\documentclass{amsart}
\usepackage{amsthm,amsmath,amsfonts}
 \newtheorem{theorem}{Theorem}[section]
 \newtheorem{corollary}[theorem]{Corollary}
 \newtheorem{lemma}[theorem]{Lemma}

 \newtheorem{defn}[theorem]{Definition}
 \newtheorem{ass}[theorem]{Assumption}
 \newtheorem{remark}[theorem]{Remark}
 \newtheorem{ex}{Example}
 \numberwithin{equation}{section}

\def\F{\mathcal{F}}
\def\R{\mathbb{R}}
\def\N{\mathbb{N}}
\def\K{\mathbf{K}}
\def\Q{\mathbf{Q}}

\def\P{\mathbf{P}}

\def\f{\mathbf{f}}
\def\z{\mathbf{z}}
\def\y{\mathbf{y}}
\def\co{{\rm co}\,}

\title[Convexity in semi-algebraic geometry]{Convexity in semi-algebraic geometry and polynomial optimization}

\author[Lasserre]{Jean B. Lasserre}

\address{%
LAAS-CNRS and Institute of Mathematics\\
University of Toulouse\\
LAAS, 7 avenue du Colonel Roche\\
31077 Toulouse C\'edex 4\\
France}

\email{lasserre@laas.fr}

\thanks{Research partially supported by the (french) ANR grant NT05-3-41612.}



\date{}
\begin{document}

\begin{abstract}
We review several (and provide new) results on the theory of 
moments, sums of squares and basic semi-algebraic sets when convexity is present. In particular, we show that under convexity, the hierarchy of semidefinite relaxations for polynomial 
optimization simplifies and has finite convergence, a highly desirable feature as convex problems are 
in principle easier to solve. In addition, if a basic semi-algebraic set $\K$ is convex but its 
defining polynomials are not, we provide two algebraic {\it certificate} of convexity
which can be checked numerically. The second is simpler and holds if a sufficient 
(and almost necessary) condition is satisfied, it also provides a new condition for $\K$ to 
have semidefinite representation.
For this we use (and extend) some of recent results from the author and Helton and Nie \cite{HN1}.
Finally, we show that when restricting to a certain class of convex polynomials, 
the celebrated Jensen's inequality in convex analysis
can be extended  to linear functionals that are not necessarily probability measures.
\end{abstract}

\maketitle

\section{Introduction}

\subsection*{Motivation} This paper is a contribution to the new emerging field of convex semi-algebraic geometry,
and its purpose is threefold: First we show that the 
moment approach for global polynomial optimization proposed in \cite{lasserre1}, 
and based on semidefinite programming (SDP), is consistent as it simplifies
and/or has better convergence properties when 
solving convex problems. 
In other words, the SDP moment approach somehow "recognizes" convexity,
a highly desirable feature for a general purpose
method because, in principle, convex problems should be easier to solve. 

We next review some recent results (and provide a new one)
on the representation of convex basic semi-algebraic sets 
by linear matrix inequalities which show how convexity permits to derive
relatively simple and {\it explicit} semidefinite representations. In doing so
we also provide a {\it certificate} of convexity for $\K$ when its defining polynomials are not convex.

Finally, we consider the important Jensen's inequality in convex analysis.
When restricting its application to a class of convex polynomials,
we provide an extension to a class of linear functionals that are not 
necessarily probability measures. 

To do so, we use (and sometimes extend) some recent results of the author
\cite{lasserre-sdr, lasserre-convex} and Helton and Nie \cite{HN1}.
We hope to convince the reader that convex semi-algebraic geometry
is indeed a very specific subarea of real algebraic geometry which
should deserve more attention from both the optimization and real algebraic
geometry research communities.

\subsection*{Background} I. Relatively recent results in the theory of moments and its dual theory of
positive polynomials have been proved useful in polynomial optimization as they provide the  
basis of a specific convergent numerical approximation scheme. Namely, one can define a hierarchy of 
semidefinite relaxations (in short SDP-relaxations) of the original optimization problem whose associated monotone sequence of optimal values converges to the global optimum. 
For a more detail account
of this approach, the interested reader is referred to e.g. Lasserre \cite{lasserre1,lasserre2}, Parrilo \cite{parrilo}, Schweighofer \cite{markus}, and the many references therein. 

Remarkably, practice seems to reveal that convergence is often fast and even finite.
However, the size of the SDP-relaxations grows rapidly with the rank in the hierarchy; 
typically the $r$-th SDP-relaxation in the hierarchy has $O(n^{2r})$ variables and semidefinite matrices
of $O(n^{r})$ sizes (where $n$ is the number of variables in the original problem).
On the other hand, it is well-known that a large class of convex optimization problems can be solved efficiently; see e.g. Ben Tal and Nemirovski \cite{bental}. 
%Hence, it is important to know whether the moment approach is also efficient in the convex case.
Therefore, as the SDP-based moment approach is dedicated to solving difficult non convex (most of 
the time NP-hard) problems, it should have the highly desirable feature 
to somehow {\it recognize} "easy" problems like convex ones. 
That is, when applied to such easy problems it should show some significant
improvement or a particular nice behavior not necessarily valid in the general case.
Notice that this is {\it not} the case of the LP-based moment-approach 
described in \cite{lasserre2,lasserre3} for which only asymptotic (and {\it not} finite) convergence occurs in general
(and especially for convex problems), a rather annoying feature.
However, for SDP-relaxations, some results of \cite{lasserre-convex} already show that indeed convexity helps as one provides specialized
representation results for convex polynomials that are nonnegative on a basic semi-algebraic set.

II. Next, in view of the 
potential of semidefinite programming techniques, 
an important issue is the characterization of convex sets that are semidefinite representable
(in short called SDr sets). A SDr set $\K\subset\R^n$ is the projection of 
a set defined by linear matrix inequalities (LMIs). That is,
\[\K:=\{x\in\R^n\::\:\exists\, y\in\R^s\:{\rm s.t.} \quad A_0+\sum_{i=1}^nx_i\,A_i+
\sum_{j=1}^sy_j\,B_j\succeq0\}\]
for some real symmetric matrices $(A_i,B_j)$ (and where $A\succeq0$ stands for $A$ is positive semidefinite).
For more details, the interested reader is referred to 
Ben Tal and Nemirovski \cite{bental}, 
Lewis et al. \cite{lewis}, Parrilo \cite{parrilo2}, 
and more recently, Chua and Tuncel \cite{chua}, 
Helton and Nie \cite{HN1,HN2}, Henrion \cite{didier} and Lasserre \cite{lasserre-sdr}.
For compact basic semi-algebraic sets %$\K\subset\R^n$ in the form
\begin{equation}
\label{setk}
\K\,:=\{x\in\R^n\::\:g_j(x)\geq0,\quad j=1,\ldots,m\:\},
\end{equation}
recent results of Helton and Nie \cite{HN1,HN2} and the author \cite{lasserre-sdr} provide sufficient conditions
on the defining polynomials $(g_j)\subset\R[X]$ for the convex hull
$\co(\K)$ ($\equiv\K$ if $\K$ is convex)  to be SDr.  
Again, an interesting issue is to analyze whether 
convexity of $\K$ (with or without concavity of the defining polynomials $(g_j)$) 
provides some additional insights and/or simplifications. 
Another interesting issue is how to detect whether a basic semi-algebraic set $\K$ is convex, or equivalently, how to obtain an algebraic {\it certificate} of convexity of $\K$ from its defining polynomials $(g_j)$.
By certificate we mean a mathematical statement that 
obviously implies convexity of $\K$, can be checked
numerically and does not require infinitely many tests. So far, and to the best of our knowledge, such a certificate does not exist.

III. The celebrated Jensen's inequality is an important result in convex analysis which
states that $E_\mu(f(x))\geq f(E_\mu(x))$ for a convex function $f:\R^n\to\R$ and a probability measure $\mu$ with $E_\mu(x)<\infty$. A third goal of this paper
is to analyze whether when restricted to a certain class of convex polynomials,
Jensen's inequality can be extended to a class of linear functionals larger than the class
of probability measures.

\subsection*{Contribution}

Concerning issue I: We first recall two previous results proved in \cite{lasserre-convex}:
 (a) the cone of convex SOS is dense (for the $l_1$-norm of coefficients) in the cone of nonnegative convex polynomials,
and (b) a convex Positivstellensatz
for convex polynomials nonnegative on $\K$ (a specialization of Putinar's Positivstellensatz).
We then analyze the role of convexity 
for the polynomial optimization problem 
\begin{equation}
\label{pbp}
\P:\quad f^*\,=\,\displaystyle\min_x\: \{\:f(x)\::\:x\in\K\:\}
\end{equation}
with $\K$ as in (\ref{setk}), and show that indeed convexity helps and makes the
SDP-relaxations more efficient. In particular, when $\K$ is convex and Slater's condition\footnote{Slater's condition holds for $\K$ in (\ref{setk}) if 
for some $x_0\in\K$, $g_j(x_0)>0$, $j=1,\ldots,m$.} holds, 
by using some recent results of Helton and Nie \cite{HN1}, we show that

(i) If the polynomials $f,(-g_j)$ are all convex and
$\nabla^2f$ is positive definite (and so $f$ is strictly convex) on $\K$,
then the hierarchy of SDP-relaxations has {\it finite} convergence.

(ii) If  $f$ and $(-g_j)$ are all 
SOS-convex (i.e. their Hessian is a SOS matrix polynomial),
then $\P$ reduces to solving a {\it single} SDP whose index in
the hierarchy is readily available.\\

Concerning II: Under certain sufficient conditions on the $(g_j)$ 
(typically some second order positive curvature conditions) 
Helton and Nie \cite{HN1,HN2} have proved that ${\rm co}\,(\K)$ (or $\K$ if convex) has
a semidefinite representation that uses Schm\"udgen or Putinar
SOS representation of polynomials positive on $\K$; see \cite{HN1,lasserre-convex}. 
Yet, in general its dimension depends on an unknown degree parameter in 
Schm\"udgen (or Putinar) SOS representation.
Our contribution is to provide a new sufficient condition for existence of a SDr
when $\K$ is compact with nonempty interior and its boundary satisfies
some nondegeneracy assumption.
It translates the geometric property of convexity of $\K$ into a 
SOS Putinar representation of some appropriate polynomial
obtained from each $g_j$. When satisfied, this representation 
provides an algebraic certificate of convexity for $\K$ and it
is almost necessary in the sense that  it always holds true
when relaxed by an arbitrary $\epsilon>0$.
It also contains as special cases Helton and Nie \cite{HN1} sufficient conditions
of SOS-convexity or strict convexity on $\partial\K$ of the $-g_j$'s, and leads to an explicit semidefinite representation of $\K$. We also provide a more general algebraic certificate based on
Stengle's Positivstellensatz, but more complex and heavy to implement
and so not very practical. In practice both certificates are
obtained by solving a semidefinite program. Therefore, because of unavoidable 
numerical inaccuracies, the certificate is valid only up to machine precision.

Concerning III, we prove that when restricting its application to the subclass of
SOS-convex polynomials, Jensen's inequality can be extended
to all linear functionals $L_\y$ (with $L_\y(1)=1$) in the 
dual cone of SOS polynomials, hence {\it not}
necessarily probability measures.

Some of the results already obtained in \cite{HN1,lasserre-sdr} and
in the present paper strongly suggest that the class of SOS-convex polynomials introduced in Helton and Nie \cite{HN1} is particularly nice and should deserve more attention.

\section{Notation, definitions and preliminary results}

Let $\R[X]$ be the ring of real polynomials in the variables $X=(X_1,\ldots,X_n)$, and let $\Sigma^2[X]\subset\R[X]$ be the subset of sums of squares (SOS) polynomials. 
Denote $\R[X]_d\subset\R[X]$ be the set of
polynomials of degree at most $d$, which forms a vector space of dimension
$s(d)={n+d\choose d}$.
If $f\in\R[X]_d$, write
$f(X)=\sum_{\alpha\in\N^n}f_\alpha X^\alpha$ in
the usual canonical basis $(X^\alpha)$, and
denote by $\f=(f_\alpha)\in\R^{s(d)}$ its vector of coefficients. Also write
$\Vert f\Vert_1 \:(=\Vert\f\Vert_1:=\sum_\alpha \vert f_\alpha\vert$) the $l_1$-norm of
$f$. Finally, denote by $\Sigma^2[X]_d\subset\Sigma^2[X]$ the subset of
SOS polynomials of degree at most $2d$. 

We use the notation 
$X$ for the variable of a polynomial $X\mapsto f(X)$ and $x$ when $x$ is a point of $\R^n$,
as for instance in $\{x\in\R^n \::\: f(x)\geq0\}$.

\subsection*{Moment matrix} With $\y=(y_\alpha)$ being a sequence indexed in the canonical basis
$(X^\alpha)$ of $\R[X]$, let $L_\y:\R[X]\to\R$ be the linear functional
\[f\quad (=\sum_{\alpha}f_{\alpha}\,X^\alpha)\quad\mapsto\quad
L_\y(f)\,=\,\sum_{\alpha}f_{\alpha}\,y_{\alpha},\]
and let $M_d(\y)$ be the symmetric matrix with rows and columns indexed in 
the canonical basis $(X^\alpha)$, and defined by:
\[M_d(\y)(\alpha,\beta)\,:=\,L_\y(X^{\alpha+\beta})\,=\,y_{\alpha+\beta},\quad\alpha,\beta\in\N^n_d\]
with $\N^n_d:=\{\alpha\in\N^n\::\:\vert \alpha\vert \:(=\sum_i\alpha_i)\leq d\}$.
\vspace{0.2cm}

\subsection*{Localizing matrix} Similarly, with $\y=(y_{\alpha})$
and $g\in\R[X]$ written
\[X\mapsto g(X)\,=\,\sum_{\gamma\in\N^n}g_{\gamma}\,X^\gamma,\]
let $M_d(g\,\y)$ be the symmetric matrix with rows and columns indexed in 
the canonical basis $(X^\alpha)$, and defined by:
\[M_d(g\,\y)(\alpha,\beta)\,:=\,L_\y\left(g(X)\,X^{\alpha+\beta}\right)\,=\,\sum_{\gamma}g_{\gamma}\,
y_{\alpha+\beta+\gamma},\]
for every $\alpha,\beta\in\N^n_d$. 
\vspace{0.2cm}

\subsection*{Putinar Positivstellensatz}
Let $Q(g)\subset\R[X]$ be the quadratic module generated by the polynomials $(g_j)\subset\R[X]$, that is,
\begin{equation}
\label{qg}
Q(g)\,:=\,\left\{\sigma_0+\sum_{j=1}^m\sigma_j\,g_j\::
\quad(\sigma_j)\subset\Sigma^2[X]\:\right\}.\end{equation}

\begin{ass}%[Putinar's assumption]
\label{assput}
$\K\subset\R^n$ is a compact basic semi-algebraic set defined as in (\ref{setk}) and
the quadratic polynomial $X\mapsto M-\Vert X\Vert^2$ belongs to $Q(g)$.
\end{ass}
\vspace{0.2cm}

Assumption \ref{assput} is not very restrictive. For instance, it holds if
every $g_j$ is affine (i.e., $\K$ is a convex polytope) or if the level
set $\{x\::\:g_j(x)\geq0\}$ is compact for some $j\in\{1,\ldots,m\}$. In addition,
if $M-\Vert x\Vert\geq0$ for all $x\in\K$, then it suffices to add the redundant quadratic constraint $M^2-\Vert x\Vert^2\ge0$ to the definition (\ref{setk}) of $\K$ and Assumption \ref{assput} will hold true.

\begin{theorem}[Putinar's Positivstellensatz \cite{putinar}]
\label{thput}
Let Assumption \ref{assput} hold.
If $f\in\R[X]$ is (strictly) positive on $\K$, then $f\in Q(g)$. That is:
\begin{equation}
\label{putinarrep}
f\,=\,\sigma_0+\sum_{j=1}^m\sigma_j\, g_j,\end{equation}
for some SOS polynomials $(\sigma_j)\subset\Sigma^2[X]$. 
\end{theorem}

\subsection{A hierarchy of semidefinite relaxations (SDP-relaxations)}

Let $\P$ be the optimization problem (\ref{pbp}) with $\K$ as in (\ref{setk})
and let $r_j=\lceil ({\rm deg}\,g_j)/2\rceil$, $j=1,\ldots,m$.
With $f\in\R[X]$ and $2r\geq \max[{\rm deg}\,f,\:\max_j2r_j]$, consider the hierarchy of semidefinite relaxations $(\Q_r)$ defined by:
\begin{equation}
\label{sdpprimal}
\Q_r:\:\quad\left\{\begin{array}{ll}
\displaystyle\inf_\y&L_\y(f)\\
\mbox{s.t.}&M_r(\y)\,\succeq0\\
&M_{r-r_j}(g_j\,\y)\,\succeq0,\qquad j=1,\ldots,m\\
&y_0\,=1
\end{array}\right.,
\end{equation}
with optimal value denoted by $\inf\Q_r$. One says that $\Q_r$ is solvable if it has an optimal solution
(in which case one writes $\inf\Q_r=\min\Q_r$).
The dual of $\Q_r$ reads
\begin{equation}
\label{sdpdual}
\Q^*_r:\:\quad\left\{\begin{array}{ll}
\displaystyle\sup&\lambda\\
\mbox{s.t.}&f-\lambda\,=\,\sigma_0+\displaystyle\sum_{j=1}^m\sigma_j\,g_j\\
&\sigma_j\in\Sigma^2[X],\quad j=0,1,\ldots,m\\
&{\rm deg}\,\sigma_0,\,{\rm deg}\,\sigma_j+{\rm deg}\,g_j\leq 2r,\quad j=1,\ldots,m
\end{array}\right.,
\end{equation}
with optimal value denoted by $\sup\Q^*_r$ (or $\max\Q^*_r$ if the $\sup$ is attained).\\

By weak duality $\sup\Q^*_r\leq\inf\Q_r$ for every $r\in\N$ and 
under Assumption \ref{assput}, 
$\inf\Q_r\uparrow f^*$ as $r\to\infty$. For a more detailed account see e.g.
 \cite{lasserre1}.

 \subsection{Convexity and SOS-convexity}
 
We first briefly recall basic facts on a multivariate convex function.
If $C\subseteq\R^n$ is a nonempty convex set, a function $f:C\to\R$ is convex on $C$ if and only if
 \[f(\lambda x+(1-\lambda)y)\leq\,\lambda f(x)+(1-\lambda)f(y),\qquad \forall\,\lambda\in (0,1),\:x,y\in C.\]
Similarly, $f$ is strictly convex on $C$ if and only if
the above inequality is strict for every $x,y\in C$, $x\neq y$, and all $\lambda\in (0,1)$.

If $C\subseteq\R^n$ is an open convex set and $f$ is twice differentiable on $C$, then $f$ is convex on $C$ if and only if
its Hessian $\nabla^2f$ is positive semidefinite on $C$ (denoted $\nabla^2f\succeq0$ on $C$).
Finally,  if $\nabla^2f$ is positive definite on $C$ (denoted $\nabla^2f\succ0$ on $C$) then $f$ is strictly convex on $C$.

\subsection*{SOS-convexity}
Helton and Nie \cite{HN1} have introduced the following interesting subclass of 
convex polynomials, called SOS-convex polynomials.\\

\begin{defn}[Helton and Nie \cite{HN1}]
\label{def1}
A polynomial $f\in\R[X]_{2d}$ is said to be SOS-convex if
$\nabla^2f$ is SOS, that is, $\nabla^2f=LL^T$ for some real matrix polynomial 
$L\in\R[X]^{n\times s}$ (for some $s\in\N$).
\end{defn}
\vspace{0.2cm}

As noted in \cite{HN1}, an important feature of SOS-convexity is that it can be can be checked numerically by solving a SDP.  They have also proved the following 
important property:\\

\begin{lemma}[Helton and Nie {\cite[Lemma 7]{HN1}}]
\label{prop}
If a symmetric matrix polynomial $P\in\R[X]^{r\times r}$ is SOS then for any $u\in\R^n$, the double integral
\[X\mapsto\quad F(X,u)\,:=\,\int_0^1\int_0^tP(u+s(X-u))\,ds\,dt\]
is also a symmetric SOS matrix polynomial in $\R[X]^{r\times r}$.
\end{lemma}
\vspace{0.2cm}

And also:
\begin{lemma}[Helton and Nie {\cite[Lemma 8]{HN1}}]
\label{prop2}
For a polynomial $f\in\R[X]$ and every $x,u\in\R^n$:
\begin{eqnarray*}
f(x)&=&f(u)+\nabla f(u)^T(x-u)\\
&&+\:(x-u)^T\underbrace{\int_0^1\int_0^t\nabla^2f(u+s(x-u))dsdt}_{F(x,u)}\,(x-u).\end{eqnarray*}
And so if $f$ is SOS-convex and $f(u)=0,\nabla f(u)=0$, then $f$ is a SOS polynomial.
\end{lemma}

\subsection{An extension of Jensen's inequality}
Recall that if $\mu$ is a probability measure on $\R^n$ with $E_\mu(x)<\infty$,
Jensen's inequality states that if $f\in L_1(\mu)$ and $f$ is convex, then
\[E_\mu(f(x))\,\geq\,f(E_\mu(x)),\]
a very useful property in many applications.

We now provide an extension of 
Jensen's inequality when one restricts its application
to the class of SOS-convex polynomials.
Namely, we may consider the linear functionals $L_\y:\R[X]_{2d}\to\R$ 
in the dual cone of $\Sigma^2[X]_{d}$, that is, vectors 
$\y=(y_\alpha)$ such that $M_d(\y)\succeq0$ and $y_0=L_\y(1)=1$;
hence $\y$ is {\it not} necessarily the (truncated) moment sequence of some probability measure $\mu$. 
Crucial in the proof is Lemma \ref{prop} of Helton and Nie.

\begin{theorem}
\label{th-jensen}
Let $f\in\R[X]_{2d}$ be SOS-convex, and let
$\y=(y_\alpha)_{\alpha\in\N^n_{2d}}$ satisfy $y_0=1$ and $M_d(\y)\succeq0$.
Then:
\begin{equation}
\label{jensen}
L_\y(f(X))\,\geq\,f(L_\y(X)),
\end{equation}
where $L_\y(X)=(L_\y(X_1),\ldots,L_\y(X_n))$.
\end{theorem}
\begin{proof}
Let $z\in\R^n$ be fixed, arbitrary, and consider the polynomial $X\mapsto f(X)-f(z)$. 
Then,
\begin{equation}
\label{jensen-1}
f(X)-f(z)\,=\,\langle \nabla f(z),X-z\rangle+\langle (X-z),F(X)(X-z)\rangle,\end{equation}
with $F:\R^n\to \R[X]^{n\times n}$ being the matrix polynomial
\[X\mapsto\quad F(X)\,:=\,\int_0^1\int_0^t\nabla^2f(z+s(X-z))\,ds\,dt.\]
As $f$ is SOS-convex, by Lemma \ref{prop}, $F$ is a SOS matrix polynomial and so the polynomial
$X\mapsto \Delta(X):=\langle (X-z),F(X)(X-z)$ is SOS, i.e., 
$\Delta\in\Sigma^2[X]$.
Then applying $L_\y$ to the polynomial $X\mapsto f(X)-f(z)$ and using (\ref{jensen-1}) yields
(recall that $y_0=1$)
\begin{eqnarray*}
L_\y(f(X))-f(z)&=&\langle\nabla f(z),L_\y(X)-z\rangle+L_\y(\Delta(X))\\
&\geq&\langle\nabla f(z),L_\y(X)-z\rangle\quad\mbox{[because $L_\y(\Delta(X))\geq0$].}
\end{eqnarray*}
As $z\in\R^n$ was arbitrary, taking $z:=L_\y(X)\,(=(L_\y(X_1),\ldots,L_\y(X_n))$ yields the desired result.
\end{proof}
\vspace{0.2cm}

As a consequence we also get:

 \begin{corollary}
Let $f$ be a convex univariate polynomial, $g\in\R[X]$ (and so 
$f\circ g\in\R[X]$). Let $d:=\lceil ({\rm deg}\,f\circ g)/2\rceil$, and let
$\y=(y_\alpha)_{\alpha\in \N^n_{2d}}$ be such that $y_0=1$ and $M_d(\y)\succeq0$. Then:
\begin{equation}
\label{cor-1}
L_\y[\,f(g(X))\,]\,\geq\,f(L_\y[\,g(X)\,]).
\end{equation}
\end{corollary}
\begin{proof}
Again let $z\in\R^n$ be fixed, arbitrary, and consider the univariate 
polynomial $Y\mapsto f(Y)-f(z)$ so that
(\ref{jensen-1}) holds. That is,
\[f(Y)-f(z)\,=\,f'(z)\,(Y-z)+ F(Y)(Y-z)^2,\]
with $F:\R\to \R[Y]$ being the univariate polynomial
\[Y\mapsto\quad F(Y)\,:=\,\int_0^1\int_0^t f"(z+s(Y-z))\,ds\,dt.\]
As $f$ is convex, $f"\geq0$, and so the univariate polynomial
$Y\mapsto F(Y)(Y-z)^2$ is nonnegative, and being univariate, is SOS.
Therefore, with $Y:=g(X)$,
\[f(g(X))-f(z)\,=\,f'(z)\,(g(X)-z)+ F(g(X))(g(X)-z)^2,\]
and so
\begin{eqnarray*}
L_\y[\,f(g(X))]-f(z)&=&f'(z)\,(L_\y[\,g(X)\,]-z)+L_\y[\,F(g(X))\,(g(X)-z)^2\,]\\
&\geq&f'(z)( L_\y[\,g(X)\,]-z)\end{eqnarray*}
and taking $z:=L_\y[g(X)]$ yields the desired result.
\end{proof}
\vspace{0.2cm}

Hence the class of SOS-convex polynomials has the very interesting property to
extend Jensen's inequality to some linear functionals that are not necessarily coming from a probability measure.

 \section{Semidefinite relaxations in the convex case}

\subsection{A convex Positivstellensatz}
Let $\K$ be as in (\ref{setk}) and define $Q_c(g)\subset \R[X]$ to be the set:
\begin{equation}
\label{qc}
Q_c(g)\,:=\,\left\{\:\sigma_0+\sum_{j=1}^m\lambda_j\,g_j\::\quad
\lambda\in\R^m_+\,;\:\sigma_0\in\Sigma^2[X],\:\sigma_0\mbox{ convex}
\:\right\}\subset Q(g).
\end{equation}
The set $Q_c(g)$ is a specialization of 
$Q(g)$ in (\ref{qg}) to the convex case, in that the weights asociated with the $g_j$'s are nonnegative scalars, i.e., SOS polynomials of degree 0, and the SOS polynomial $\sigma_0$ is convex. 
In particular, every $f\in Q_c(g)$ is nonnegative on $\K$.
Let $\F_\K\subset\R[X]$ be the convex cone
of convex polynomials nonnegative on $\K$.\\

\begin{theorem}[Lasserre \cite{lasserre-convex}]
\label{thmain2}
Let $\K$ be as in (\ref{setk}), Slater's condition hold
and $g_j$ be concave for every $j=1,\ldots,m$.

Then with $Q_c(g)$ as in (\ref{qc}), the set $Q_c(g)\cap \F_\K$ is dense in $\F_\K$ for the $l_1$-norm $\Vert \cdot\Vert_1$. In particular, 
if $\K=\R^n$ (so that $\F_{\R^n}=:\F$ is now the set of nonnegative convex polynomials), then
$\Sigma^2[X]\cap \F$ is dense in $\F$.
%; that is, almost all nonnegative convex polynomials are SOS.
\end{theorem}
\vspace{0.2cm}

Theorem \ref{thmain2} states that if $f$ is convex and nonnegative on $\K$
(including the case $\K\equiv \R^n$)  then one may approximate $f$ by a sequence
$\{f_{\epsilon r}\}\subset Q_c(g)\cap \F_\K$
with $\Vert f-f_{\epsilon r}\Vert_1\to 0$ as $\epsilon\to 0$
(and $r\to\infty$). For instance, with $r_0:=\lfloor ({\rm deg}\,f)/2\rfloor +1$,
\begin{eqnarray}
\nonumber
X\mapsto f_{\epsilon r}(X)&:=&f+\epsilon (\theta_{r_0}(X)+\theta_{r}(X)),\quad\mbox{ with }\\
\label{2}
X\mapsto \theta_r(X)&:=&1+\sum_{k=1}^r \sum_{i=1}^n \frac{X_i^{2k}}{k{\rm !}}
\qquad r\geq r_\epsilon,\end{eqnarray}
for some $r_\epsilon$; see Lasserre \cite{lasserre-convex} for details.
Observe that Theorem \ref{thmain2} provides $f$ with a {\it certificate}
of nonnegativity on $\K$.
Indeed, let $x\in\K$ be fixed arbitrary. Then as $f_{\epsilon r}\in Q_c(g)$
one has $f_{\epsilon r}(x)\geq0$. Letting $\epsilon\downarrow 0$ yields
$0\leq \lim_{\epsilon\to 0}f_{\epsilon r}(x)=f(x)$. And as $x\in\K$ was arbitray,
$f\geq0$ on $\K$.

Theorem \ref{thmain2} is a convex (weak) version of Theorem \ref{thput} (Putinar's Positivstellensatz) where one replaces the quadratic module $Q(g)$ with its subset
$Q_c(g)$. We call  it a {\it weak} version of Theorem \ref{thput} because
it invokes a density result (i.e. $f_{\epsilon r}\in Q_c(g)$ whereas $f$ might not be an element of $Q_c(g)$).
Notice that $f$ is allowed to be nonnegative (instead of strictly positive) on $\K$ and
$\K$ need {\it not} be compact; recall that extending 
Theorem \ref{thput} to non compact basic semi-algebraic sets $\K$ and to polynomials $f$ nonnegative on $\K$
is hopeless in general; see Scheiderer \cite{claus}.\\

\begin{corollary}
\label{special}
Let $\K$ be as in (\ref{setk}), $f\in\R[X]$ with $f^*:=\min_x \{f(x)\::\:x\in\K\}$ and let 
$d:=\max[\lceil ({\rm deg}\,f)/2\rceil,\max_j\lceil ({\rm deg}\,g_j)/2\rceil\,]$.
Consider the simplified SDP-relaxation
 \begin{equation}
\label{newsdpprimal}
\widehat{\Q}:\:\quad\left\{\begin{array}{ll}
\displaystyle\inf_\y&L_\y(f)\\
\mbox{s.t.}&M_d(\y)\,\succeq0\\
&L_\y(g_j)\,\geq0,\qquad j=1,\ldots,m\\
&y_0\,=1
\end{array}\right.
\end{equation}
and its dual
\begin{equation}
\label{newsdpdual}
\widehat{\Q}^*:\:\quad\left
\{\begin{array}{ll}
\displaystyle\sup_{\gamma,\sigma_0,\lambda}&\gamma\\
\mbox{s.t.}&f-\gamma\,=\,\sigma_0+\displaystyle\sum_{j=1}^m\lambda_j\,g_j\\
&\sigma_0\in\Sigma^2[X]_d;\:\lambda_j\geq0,\quad j=1,\ldots,m
\end{array}\right.
\end{equation}
\indent
{\rm (a)} If  $f-f^*\in Q_c(g)$
%, i.e., if $f-f^*\,=\,\sigma_0+\sum_{j=1}^m\lambda_j\,g_j$ with $\sigma_0\in\Sigma^2[X]$,
% and $\lambda\in\R^m_+$,  
then the SDP-relaxation 
 $\widehat{\Q}$ and its dual $\widehat{\Q}^*$ are exact.\\
 
 {\rm (b)} If $f,-g_j\in\R[X]$ are convex, $j=1,\ldots,m$, and if $\y$ is an optimal solution of $\widehat{\Q}$
 which satisfies
 \begin{equation}
 \label{rank}
  {\rm rank}\,M_d(\y)\,=\, {\rm rank}\,M_{d-1}(\y),\end{equation}
 then $\widehat{\Q}$ is exact and 
 $x^*:=(L_\y(X_i))\in\K$ is a (global) minimizer of $f$ on $\K$.
 \end{corollary}
\begin{proof}
(a) If $f-f^*\in Q_c(g)$, i.e., if $f-f^*=\sigma_0+\sum_{j=1}^m\lambda_jg_j$, 
with $\sigma_0\in\Sigma^2[X]_d$ and $\lambda\in\R^m_+$, 
the triplet $(f^*,\sigma_0,\lambda)$ is a feasible solution of $\widehat{\Q}^*$
with value $f^*$. Therefore, as $\sup\widehat{\Q}^*\leq\inf\widehat{\Q}\leq f^*$, the SDP-relaxation $\widehat{\Q}$ and its dual $\widehat{\Q}^*$ are exact. In fact, $(f^*,\sigma_0,\lambda)$ is an optimal solution of $\widehat{\Q}^*$.

(b) If $\y$ satisfies the rank condition (\ref{rank}) then 
by the {\it flat extension} theorem of Curto and Fialkow \cite{curto},
$\y$ is the (truncated) moment sequence of an atomic probability measure $\mu$ on $\R^n$,
say $\mu=\sum_{k=1}^s\lambda_k\delta_{x(k)}$ with 
$s={\rm rank}\,M_d(\y)$, $0<\lambda_k\leq 1$, $\sum_k\lambda_k=1$,
and $\delta_{x(k)}$ being the Dirac measure at $x(k)\in\R^n$, $k=1,\ldots,s$.
Let $x^*:=\sum_k\lambda_kx(k)=(L_\y(X_i))\in\R^n$.
Then $f^*\geq L_\y(f)$ and by convexity of $f$, $L_\y(f)=\sum_k\lambda_kf(x(k))\geq f(\sum_k\lambda_kx(k))=f(x^*)$.
Similarly, by convexity of $-g_j$,
$0\leq L_\y(g_j)=\sum_k\lambda_kg_j(x(k))\leq g_j(\sum_k\lambda_kx(k))=g_j(x^*)$, $j=1,\ldots,m$.
Therefore, $x^*\in\K$ and as $f(x^*)\leq f^*$, $x^*$ is a global minimizer of $f$ on $\K$.
\end{proof}
\vspace{0.2cm}

Notice that $\K$ in Corollary \ref{special} need not be compact. Also,
Corollary \ref{special}(b) has practical value because in general 
one does not know whether $f-f^*\in Q_c(g)$
(despite that in the convex case, $f-f^*\in\mathcal{F}_\K$ and $Q_c(g)\cap\mathcal{F}_\K$
is dense in $\mathcal{F}_\K$). However, one may still solve $\widehat{\Q}$ and check
whether the rank condition (\ref{rank}) is satisfied.
If in solving $\widehat{\Q}_r$, the rank condition (\ref{rank}) is not satisfied, then
other sufficient conditions can be exploited as we next see.

\subsection{The SOS-convex case}
Part (a) of the following result is already contained in 
Lasserre \cite[Cor. 2.5]{lasserre-convex}.\\

\begin{theorem}
\label{thm2}
Let $\K$ be as in (\ref{setk}) and Slater's condition hold.
Let $f\in\R[X]$ be such that
$f^*:=\inf_x\{f(x)\::\:x\in\K\}=f(x^*)$ for some $x^*\in\K$.
If $f$ is SOS-convex and $-g_j$ is SOS-convex
for every $j=1,\ldots,m$, then:

{\rm (a)} $f-f^*\in Q_c(g)$.

{\rm (b)} The simplified SDP-relaxation $\widehat{\Q}$ 
in (\ref{newsdpprimal}) and its dual (\ref{newsdpdual}) are exact and solvable.
If $\y$ is an optimal solution of $\widehat{\Q}$ then
$x^*:=(L_\y(X_i))\in\K$ is a global minimizer of $f$ on $\K$.
\end{theorem}
\begin{proof}
(a) is proved in \cite[Cor. 2.5]{lasserre-convex}.
(b) That  $\widehat{\Q}$ is exact follows from (a) and Corollary \ref{special}(a). 
Hence it is solvable (e.g. take $\y$ to be the moment sequence associated with 
the Dirac measure at a global minimizer $x^*\in\K$). So let
$\y$ be an optimal solution of $\widehat{\Q}$, hence with $f^*=L_\y(f)$. As
$-g_j$ is SOS-convex for every $j$, then by Theorem \ref{th-jensen},
$0\leq L_\y(g_j)\leq g_j(x^*)$ with $x^*:=(L_\y(X_i))$ and so $x^*\in\K$.
Similarly, as $f$ is SOS-convex, we also have $f^*=L_\y(f)\geq f(x^*)$ which proves that 
$f(x^*)=f^*$ and $x^*$ is a global minimizer of $f$ on $\K$.
Finally, as by (a) $f-f^*\in Q_c(g)$ then $\widehat{\Q}^*$ is exact and solvable.
\end{proof}
\vspace{0.2cm}

(Again notice that $\K$ in Theorem \ref{thm2} need not be compact.)
So the class of SOS-convex polynomials is particularly interesting.
Not only Jensen's inequality can be extended to some linear functionals that are
not coming from a probability measure, but one may also solve SOS-convex optimization problems
$\P$ in (\ref{pbp}) (i.e. with $f$ and $\K$ defined with SOS-convex polynomials) by
solving the single semidefinite program (\ref{newsdpprimal}). 

Notice that a self-concordant\footnote{The self-concordance property
introduced in \cite{nesterov} is fundamental in the design and efficiency
of interior point methods for convex programming.} logarithmic barrier function exists for (\ref{newsdpprimal})
whereas
the logarithmic barrier function with barrier parameter $\mu$: 
\begin{equation}
\label{log}
x\mapsto \phi_\mu(x)\,:=\,\mu\,f(x)-\sum_{j=1}^m\ln\, (-g_j(x)),\end{equation}
associated with $\P$, is not self-concordant in general.
Therefore, despite (\ref{newsdpprimal}) involves additional variables (a lifting), solving (\ref{newsdpprimal}) via an
interior point method might be more efficient than solving $\P$
by using the logarithmic barrier function (\ref{log}) with no lifting.
In addition, all SOS-convex polynomials nonnegative on $\K$ and which attain
their minimum on $\K$, belong to $Q_c(g)$, a very specific 
version of Putinar Positivstellensatz (as $f$ is only nonnegative and 
$\K$ need not be compact).

\subsection{The strictly convex case}

If $f$ or some of the $-g_j$'s is not SOS-convex 
but $\nabla^2f\succ0$ (so that $f$ is strictly convex)
and $-g_j$ is convex for every $j=1,\ldots,m$, then inspired
by a nice argument from Helton and Nie \cite{HN1} for existence of a semidefinite representation of convex sets, one obtains the following result.\\

\begin{theorem}
Let $\K$ be as in (\ref{setk}) and let Assumption \ref{assput} and Slater's condition hold.
Assume that $f,-g_j\in\R[X]$ are convex, $j=1,\ldots,m$, with $\nabla^2f\succ0$ on $\K$.
%and $-\nabla^2g_j(x)\succeq0$, $j=1,\ldots,m$.

Then the hierarchy of SDP-relaxations defined in 
(\ref{sdpprimal}) has finite convergence. That is,
$f^*=\sup\Q^*_r\,=\,\inf\Q_r$ for some index $r$. In addition,
$\Q_r$ and $\Q^*_r$ are solvable so that $f^*=\max\Q^*=\min\Q_r$.
\end{theorem}
\begin{proof}
Let $x^*\in\K$ be a global minimizer (i.e. $f^*=f(x^*)$).
As Slater's condition holds, there exists a vector of Karush-Kuhn-Tucker (KKT) multipliers
$\lambda\in\R^m_+$ such that the (convex) Lagrangian $L_f\in\R[X]$ defined by
\begin{equation}
\label{lagrangian}
X\mapsto L_f(X)\,:=\,f(X)-f^*-\sum_{j=1}^m\lambda_j\,g_j(X)\end{equation}
has a global minimum at $x^*\in\K$, i.e., $\nabla L_f(x^*)=0$. In addition,
$\lambda_jg_j(x^*)=0$ for every $j=1,\ldots,m$ and $L_f(x^*)=0$.  Then,
by Lemma \ref{prop2},
\[L_f(X)\,=\,\left\langle (X-x^*),F(X,x^*)(X-x^*)\right\rangle\]
with 
\[F(X,x^*)\,:=\,\left(\int_0^1\int_0^t\nabla^2L_f(x^*+s(X-x^*))\,ds\,dt\right).\]
Next, let $I_n$ be the $n\times n$ identity matrix.
As $\nabla^2f\succ0$ on $\K$, continuity of the (strictly positive) 
smallest eigenvalue of $\nabla^2f$
and compactness of $\K$ yield that $\nabla^2f \succeq\delta I_n$ on $\K$,
for some $\delta>0$. Next, as $-g_j$ is convex for every $j$, and in view of the definition (\ref{lagrangian})
of $L_f$, $\nabla^2 L_f\succeq\nabla^2f\succeq\delta I_n$ on $\K$. 
Hence for every $\xi\in\R^n$, $\xi^TF(x,x^*)\xi\geq\delta \int_0^1\int_0^t\xi^T\xi dsdt=\frac{\delta}{2} \xi^T\xi$, and so
$F(x,x^*) \succeq\frac{\delta}{2}\,I_n$ for every $x\in\K$.
Therefore,
by the matrix polynomial version of Putinar Positivstellensatz,
\[F(X,x^*)\,=\,F_0(X)+\sum_{j=1}^mF_j(X)\,g_j(X),\]
for some real SOS matrix polynomials 
 $X\mapsto F_j(X)=L_j(X)L_j(X)^T$ (for some apppropriate
$L_j\in\R[X]^{n\times p_j}$), $j=0,\ldots,m$. See
Helton and Nie \cite{HN1}, Kojima and Maramatsu \cite{kojima}, Hol and Scherer \cite{hol}.
But then
\[X\mapsto \left\langle (X-x^*),F_j(X,x^*)(X-x^*)\right\rangle\,=\,\sigma_j(X)\in\Sigma^2[X],\qquad j=0,\ldots,m\]
and so
\begin{eqnarray*}
f(X)-f^*&=&L_f(X)+\sum_{j=1}^m\lambda_jg_j(X)\\
&=&\sigma_0(X)+\sum_{j=1}^m(\lambda_j+\sigma_j(X))\,g_j(X).
\end{eqnarray*}
Let $2s$ be the maximum degree of the SOS polynomials $(\sigma_j)$.
Then $(f^*,\{\sigma_j+\lambda_j\})$ is a feasible solution of the SDP-relaxation
$\Q^*_r$  in (\ref{sdpdual}) with $r:=s+\max_jr_j$. Therefore, as $\sup\Q^*_r\leq \inf\Q_r\leq f^*$,
the SDP-relaxations $\Q_r$ and $\Q^*_r$ are exact, finite convergence occurs and $\Q^*_r$ is solvable.
But this also implies that $\Q_r$ is solvable (take $\y$ to be the moment sequence of the 
Dirac measure $\delta_{x^*}$ at any global minimizer $x^*\in\K$).
\end{proof}
\vspace{0.2cm}

When compared to Theorem \ref{thm2} for the SOS-convex case, in the strictly convex
case the simplified SDP-relaxation $\widehat{\Q}$ in (\ref{newsdpprimal}) is not guaranteed
to be exact. However, finite convergence still occurs for the SDP-relaxations
($\Q_r$) in (\ref{sdpprimal}).

\begin{remark}
{\rm It is worth emphasizing that in general, the hierarchy of LP-relaxations 
(as opposed to SDP-relaxations) defined in 
\cite{lasserre3} and based on Krivine's representation \cite{krivine,vasilescu} for polynomials positive on $\K$,
{\it cannot} have finite convergence, especially in the convex case! For more details, the interested reader is referred to 
\cite{lasserre2,lasserre3}. Therefore, and despite LP software packages can
solve LP problems of very large size, using LP-relaxations 
does not seem a good idea even  for solving a convex polynomial optimization problem.} 
\end{remark}

\section{Convexity and semidefinite representation of convex sets}
We now consider the semidefinite representation of 
convex sets. First recall the following result.
\begin{theorem}[Lasserre \cite{lasserre-sdr}]
\label{lagrangiansos}
Let $\K$ in (\ref{setk}) be compact with
$g_j$ concave, $j=1,\ldots,m$, and assume that Slater's condition holds. 
If the Lagrangian polynomial
$L_f$ in (\ref{lagrangian}) associated  with every {\it linear}
polynomial $f\in\R[X]$ is SOS, then with $d:=\max_j\lceil ({\rm deg}\,g_j)/2\rceil$, the set
\begin{equation}
\label{lfsos}
\Omega\,:=\,\left\{(x,\y)\in\R^n\times \R^{s(2d)}\::\:\left\{\begin{array}{ll}
M_d(\y)&\succeq0\\
L_\y(g_j)&\geq 0,\quad j=1,\ldots,m\\
L_\y(X_i)&=x_i,\quad i=1,\ldots,n\\
y_0&=1\end{array}\right.\right.\end{equation}
is a semidefinite representation of $\K$. 
\end{theorem}
\vspace{0.2cm}

Next, Helton and Nie \cite{HN1,HN2} have provided several interesting 
second-order positive curvature (sufficient and necessary)  conditions on the defining polynomials
$(g_j)$ for $\K$ (or its convex hull ${\rm co}\,(\K)$) to have a SDr.
In particular (recall that $r_j=\lceil ({\rm deg}\,g_j)/2\rceil$ for every $j=1,\ldots,m$):\\

\begin{theorem}[Helton and Nie \cite{HN1}]
\label{suffhn1}
Let $\K$ in (\ref{setk}) be convex, Asssumption \ref{assput} hold, and assume that Slater's condition holds and $g_j$ is concave on $\K$, $j=1,\ldots,m$.

{\rm (a)} If $-g_j$ is SOS-convex for every $j=1,\ldots,m$, then for every linear $f\in\R[X]$, 
the associated Lagrangian $L_f$ 
(\ref{lagrangian}) is SOS and the set $\Omega$ in
(\ref{lfsos}) is a semidefinite representation of $\K$.

{\rm (b)} If every $-g_j$ is either SOS-convex or satisfies
$-\nabla^2g_i\succ0$ on $\K\cap\{x\::\:g_j(x)=0\}$, then there exists $r\in\N$
such that the set
\begin{equation}
\label{thmain-2}
\Omega\,:=\,\left\{(x,\y)\in\R^n\times \R^{s(2r)}\::\:\left\{\begin{array}{ll}
M_r(\y)&\succeq0\\
M_{r-r_j}(g_j\,\y)&\succeq 0,\quad j=1,\ldots,m\\
L_\y(X_i)&=x_i,\quad i=1,\ldots,n\\
y_0&=1\end{array}\right.\right\}\end{equation}
is a semidefinite representation of $\K$. 
\end{theorem}
\vspace{0.2cm}

See \cite[Theor. 6, and 9]{HN1}. This follows from the fact
that the Hessian $\nabla^2L_f$ associated with a linear $f\in\R[X]$ 
has a Putinar representation in terms of SOS matrix polynomials, and with degree 
of the weights bounded uniformly in $f$. In principle, the degree parameter $d$ in Theorem
\ref{suffhn1}(b) may be computed by solving a hierarchy of semidefinite programs.
Some other (more technical) weaker second-order positive curvature sufficient conditions 
(merely for existence of a SDr) are also 
provided in  \cite{HN1,HN2} but the semidefinite representation is not explicit
any more in terms of the defining polynomials $(g_j)$. Notice that if $\K$ is compact but Assumption \ref{assput} does not hold, then one still obtains a semidefinite representation for $\K$ but more complicated as it is now based on 
Schm\"udgen's representation \cite{schmudgen} instead of Putinar's representation; see \cite[Theor. 5]{HN1}.

We next provide a sufficient condition in the case where $\K$ is convex but its
defining polynomials  $(-g_j)$ are {\it not} necessarily convex. Among its distinguishing features,
it is checkable numerically, contains Theorem \ref{suffhn1} as a special case and leads to
the explicit semidefinite representation (\ref{thmain-2}) of $\K$.

\subsection{Algebraic certificate of convexity}
We first present the following characterization of convexity when
$\K$ is closed, satisfies a nondegeneracy assumption on its boundary, and Slater's condition holds.\\

\begin{lemma}
\label{lemmaconvex}
Let $\K$ be as in (\ref{setk}) (hence closed), Slater's condition hold and assume that for every
$j=1,\ldots,m$, $\nabla g_j(y)\neq0$ if $y\in\K$ and $g_j(y)=0$.
Then $\K$ is convex if and only if for every $j=1,\ldots,m$,
\begin{equation}
\label{statconvex}
\langle\nabla g_j(y),x-y\rangle \geq0,\qquad \forall x\,\in\K\mbox{ and }\forall\,y\in\K\mbox{ with }g_j(y)\,=\,0.\end{equation}
\end{lemma}
\begin{proof}
The {\it only if part} is obvious. Indeed if 
$\langle\nabla g_j(y),x-y\rangle <0$ for some $x\in\K$ and $y\in\K$ with
$g_j(y)=0$, then there is some $\overline{t}>0$ such that $g_j(y+t(x-y))<0$ for all $t\in (0,\overline{t})$ and so the point
$x':=tx+(1-t)y$ does not belong to $\K$, which in turn implies that $\K$ is not convex.

For the {\it if part}, (\ref{statconvex}) implies that at every point of the boundary, there exists a supporting hyperplane for $\K$. As $\K$ is closed with nonempty interior,
the result follows from \cite[Theor. 1.3.3]{schneider}\footnote{The author is grateful to L. Tuncel for 
providing us with the reference \cite{schneider}.}.
\end{proof}
\vspace{0.2cm}

The nondegeneracy assumption is crucial as demonstrated in the following simple example kindly provided by an anonymous referee:
\begin{ex}
{\small{\rm
\label{ex0}Consider the non convex set $\K\subset\R^2$ defined by:
\[\K\,:=\,\{\,x\in\R^2\::\:(1-x_1^2+x_2^2)^3\geq0,\:10-x_1^2-x_2^2\geq0\:\}\]
Then it is straightforward to see
that (\ref{statconvex}) is satisfied. This is because 
$\nabla g_1$ vanishes on the piece of boundary  determined by
$g_1(x)=0$.}}
\end{ex}
\vspace{0.2cm}

Next, using the above characterization (\ref{statconvex}), we provide an algebraic certificate of convexity.\\

\begin{corollary}[Algebraic certificate of convexity]
\label{iff}
Let $\K$ be as in (\ref{setk}), Slater's condition hold and assume that for every
$j=1,\ldots,m$, $\nabla g_j(y)\neq0$ if $y\in\K$ and $g_j(y)=0$.
Then $\K$ is convex if and only if for every $j=1,\ldots,m$,
\begin{equation}
\label{certif-alg}
h_j(X,Y) \langle \nabla g_j(Y),X-Y\rangle\,=\,\langle\nabla g_j(Y),X-Y\rangle^{2l}+\theta_j(X,Y)+\varphi_j(X,Y) g_j(Y),
\end{equation}
for some integer $l\in\N$,  some polynomial $\varphi_j\in\R[X,Y]$ and some
polynomials $h_j,\theta_j$ in the preordering\footnote{
The preordering of $\R[X]$ generated by a family $(g_1,\ldots,g_m)\subset\R[X]$ is the set of polynomials $\{p\::\:p=\sum_{J\subseteq\{1,\ldots,m\}}\sigma_J(\prod_{j\in J}g_j),\:\mbox{with }\sigma_J\in\Sigma^2[X]\}$.} of $\R[X,Y]$
generated by the family of polynomials $(g_k(X),g_p(Y))$, $k,p\in\{1,\ldots,m\}$, $p\neq j$.
\end{corollary}
\begin{proof}
By Lemma \ref{lemmaconvex}, $\K$ is convex if and only if for every $j=1,\ldots ,m$,
the polynomial $(X,Y)\mapsto \langle\nabla g_j(Y),X-Y\rangle$ is nonnegative on the set 
$\Omega_j$ defined by:
\begin{equation}
\label{omegaj}
\Omega_j\,:=\,\{(x,y)\in\K\times\K\::\:g_j(y)=0\:\}.\end{equation}
Equivalently, $\K$ is convex if and only if for every $j=1,\ldots,m$:
\begin{eqnarray*}
\emptyset\,=\,\left\{(x,y)\in\R^n\::\: \right.&&(x,y)\in\K\times\K\,;\quad g_j(y)\,=\,0\,;\\
&&\left.\langle\nabla g_j(y),x-y\rangle\leq0\,;\:\langle\nabla g_j(y),x-y\rangle\neq0\right\}.\end{eqnarray*}
Then (\ref{certif-alg}) follows from Stengle's Positivstellensatz \cite[Theor. 4.4.2, p. 92]{roy}.
\end{proof}
\vspace{0.2cm}

Observe that Corollary \ref{iff} provides an algebraic certificate of convexity when $\K$
is closed with nonempty interior and a nondegeneracy assumption holds on its boundary.
If one fixes an a priory bound $s$ on $l\in\N$ and on the degree of $h_j,\theta_j$ and $\varphi_j$, then checking whether
(\ref{certif-alg}) holds reduces to solving a semidefinite program. If $\K$ is convex,
by increasing $s$, eventually one would obtain such a certificate if
one could solve semidefinite programs exactly. In practice, and because of unavoidable numerical inaccuracies, one only obtains 
a numerical approximation of the optimal value and so, a certificate valid {\it up to machine precision} only.

However, implementing such a procedure is extremely costly because 
one has potentially $2\times 2^{m}$ unknown SOS polynomials to define
$h_j$ and $\theta_j$ in (\ref{certif-alg})! Therefore,
it is highly desirable to provide a less costly certificate but with no guarantee to hold 
for every $\K$ as in Corollary \ref{iff}.

In particular one only considers compact sets $\K$. Indeed,
if $\K$ is compact, one has the following result (recall that $g_0\equiv 1$).\\

\begin{lemma}
\label{lemmageom}
Let $\K$ be convex, Assumption \ref{assput} and Slater's condition hold. 
Assume that for every $j=1,\ldots,m$,
$\nabla g_j(y)\neq0$ if $y\in\K$ and $g_j(y)=0$.
Then for every $\epsilon>0$ and every $j=1,\ldots,m$:
\begin{eqnarray}
\nonumber
\left\langle\, \nabla g_j(Y),X-Y\,\right\rangle+\epsilon &=&
\sum_{k=0}^m\sigma_{jk}(X,Y)\, g_k(X)
+\sum_{k=0,k\neq j}^m\psi_{jk}(X,Y) \,g_k(Y)\\
\label{everyj}
&&+\,\psi_j(X,Y)\,g_j(Y),
\end{eqnarray}
for some SOS polynomials $(\sigma_{jk})$ and $(\psi_{jk})_{k\neq j}\subset\Sigma^2[X,Y]$, and some polynomial $\psi_j\in\R[X,Y]$. 
\end{lemma}
\begin{proof}
By Lemma \ref{lemmaconvex}, for every $j=1,\ldots,m$,
and every $x,y\in\K$ such that $g_j(y)=0$, (\ref{statconvex}) holds
%\[\langle\nabla g_j(y),x-y\rangle \geq0,\qquad \forall x\in\K,\]
and therefore, for every $j=1,\ldots,m$, 
\begin{equation}
\label{pos}
\langle \nabla g_j(y),x-y\rangle\,+\,\epsilon \,>\,0\qquad \forall (x,y)\in\Omega_j,
\end{equation}
where $\Omega_j$ has been defined in (\ref{omegaj}).
As $\K$ satisfies Assumption \ref{assput} then so does 
$\Omega_j$ for every $j=1,\ldots,m$. Hence (\ref{everyj}) follows from
(\ref{pos}) and Theorem \ref{thput}.
\end{proof}
\vspace{0.2cm}

Therefore, inspired by Lemma \ref{lemmageom}, introduce the following condition:
\begin{ass}[Certificate of convexity]
\label{geom}
For every $j=1,\ldots,m$, (\ref{everyj}) holds with $\epsilon=0$. Then
let $d_j\in\N$ be such that $2d_j$ is larger than the maximum degree of the polynomials
$\sigma_{jk}g_k,\psi_{jk}g_k,\psi_jg_j\in\R[X,Y]$ in (\ref{everyj}), $j=1,\ldots,m$.
\end{ass}
\vspace{0.2cm}

When $\K$ is closed (and not necessarily compact), Slater's condition holds and the nondegeneracy assumption
on the boundary holds (i.e., $\nabla g_j(y)\neq0$ if $y\in\K$ and $g_j(y)=0$)
Assumption \ref{geom} is indeed a certificate of convexity because
then (\ref{statconvex}) holds for every $x,y\in\K$ with $g_j(y)=0$, and by Lemma \ref{lemmaconvex},
$\K$ is convex. It translates the 
geometric property of convexity of $\K$ into an algebraic SOS Putinar representation of
the polynomial $(X,Y)\mapsto \langle\nabla g_j(Y),X-Y\rangle$ nonnegative on $\Omega_j$, $j=1,\ldots,m$. 
On the other hand, if $\K$ is convex and Assumption \ref{assput}, Slater's condition and the nondegeneracy assumption all hold, then Assumption \ref{geom} is almost necessary as, by Lemma \ref{lemmageom},
(\ref{everyj}) holds with $\epsilon>0$ arbitrary. 

With $d_j$ fixed a priori, checking whether (\ref{everyj}) hold with
$\epsilon=0$ can be done numerically. (However, again it provides a certificate of convexity
valid {\it up to machine precision} only.)
For instance, for every $j=1,\ldots,m$, it suffices to solve the semidefinite program
(recall that $r_k=\lceil ({\rm deg}\,g_k)/2\rceil$, $k=1\ldots,m$)
\begin{equation}
\label{test}
\left\{\begin{array}{ll}
\rho_j:=&\displaystyle\min_{\z}\:L_\z(\langle\nabla g_j(Y),X-Y\rangle)\\
\mbox{s.t.}&M_{d_j}(\z)\succeq0\\
&M_{d_j-r_k}(g_k(X)\,\z)\succeq0,\quad k=1,\ldots,m\\
&M_{d_j-r_k}(g_k(Y)\,\z)\succeq0,\quad k=1,\ldots,m;\,k\neq j\\
&M_{d_j-r_j}(g_j(Y)\,\z)\,=\,0\\
&y_0=1
\end{array}\right..\end{equation}
If $\rho_j=0$ for every $j=1,\ldots,m$, then 
Assumption \ref{geom} holds.
This is in contrast to the PP-BDR 
property in \cite{lasserre-convex} that cannot be checked numerically as it involves infinitely many linear polynomials $f$.

\begin{remark}
\label{numeric}
{\rm Observe that the usual rank condition (\ref{rank}) used as a 
stopping criterion to detect whether (\ref{test}) is exact (i.e. $\rho_1=0$), cannot be satisfied 
in solving (\ref{test}) with primal dual interior point methods (as in the SDP-solvers
used by GloptiPoly) because one tries to find an optimal solution $\z^*$ in the {\it relative interior}
of the feasible set of (\ref{test}) and this gives maximum rank to the moment matrix $M_{d_j}(\z^*)$. 
Therefore, in the context of (\ref{test}), if indeed $\rho_j=0$ then $\z^*$ corresponds to the moment vector
of some probability measure $\mu$ supported on the set of points $(x,x)\in \K\times\K$ that satisfy
$g_j(x)=0$ (as indeed $L_{\z^*}(\langle\nabla g_j(Y),X-Y)\rangle)=0=\rho_j$). Therefore $\rho_j=0$ as $d_j$ increases but the rank of $M_{d_j}(\z^*)$ does not stabilize because
$\mu$ is not finitely supported. In particular, a good candidate $\z^*$ for optimal solution is
the moment vector of the probability measure uniformly distributed on
the set $\{(x,x)\in\K\times\K\::\:g_j(x)=0\}$.

Alternatively, if $\rho_j\approx0$ and the dual of (\ref{test}) has an optimal solution
$(\sigma_{jk},\psi_{jk},\psi_j)$, then in some cases
one may check if (\ref{everyj}) holds exactly after appropriate rounding of coefficients of the solution.
But in general, obtaining an exact certificate (i.e., 
$\rho_j=0$ in the primal or (\ref{everyj}) with $\epsilon=0$ in the dual)
numerically is hopeless.
}
\end{remark}

\begin{ex}
\label{ex1}
{\rm {\small Consider the following simple illustrative example in $\R^2$:
\begin{equation}
\label{setexample}
\K\,:=\,\{\,x\in\R^2\::\: x_1x_2-1/4\geq0;\: 0.5-(x_1-0.5)^2-(x_2-0.5)^2\geq0\,\}\end{equation}
Obviously $\K$ is convex but its defining polynomial $x\mapsto g_1(x):=x_1x_2-1/4$ is not concave
whereas $x\mapsto g_2(x):=0.5-(x_1-0.5)^2-(x_2-0.5)^2$ is. 

With $d_1=3$, solving (\ref{test}) using GloptiPoly 3\footnote{GloptiPoly 3 (a Matlab based public software)
is an extension of GloptiPoly \cite{acm} to solve the generalized problem of moments described in \cite{rio}.
For more details see {\tt www.laas.fr/$\sim$henrion/software/}.} 
yields the optimal value 
$\rho_1\approx-4.58.10^{-11}$ which, in view of the machine precision for the SDP solvers used in GloptiPoly, could be considered to be zero, but of course with no guarantee. However, and according to Remark \ref{numeric},
we could check that 
(again up to machine precision) for every $\alpha\in\N^n$ with $\vert\alpha\vert\leq 2d_j$,
$z^*_{\alpha,\alpha}=z^*_{2\alpha,0}$ and $z^*_{\alpha,0}=z^*_{0,\alpha}$.
In addition, because of symmetry, $z_{\alpha,\beta}=z_{\alpha' ,\beta'}$
whenever $\alpha'_1=\alpha_2$ and $\alpha'_2=\alpha_1$ (and similarly for $\beta$ and $\beta'$).
Indeed for moments of order $1$ we have
$z^*_{\alpha,\beta}=(0.5707,0.5707,0.5707,0.5707)$ and for moments of order $2$,
\[z^*_{\alpha,\beta}=(0.4090,0.25,0.4090,0.25,0.4090, 0.25, 0.4090, 0.4090, 0.25,0.4090).\]
For $j=2$ there is no test to perform because $-g_2$ being quadratic and convex yields
\begin{equation}
\label{newtest}
\langle\nabla g_2(Y),X-Y\rangle \,=\,g_2(X)-g_2(Y)+\underbrace{(X-Y)^T(-\nabla^2g_2(Y))(X-Y)}_{SOS}
\end{equation}
which is in the form (\ref{everyj}) with $d_2=1$.
}}
\end{ex}
\vspace{0.2cm}

We next show the role of Assumption \ref{geom} in obtaining a semidefinite representation of $\K$.\\

\begin{theorem}
\label{thmain}
Let Assumption \ref{assput} and Slater's condition hold.
Moreover, assume that for every $j=1,\ldots,m$, 
$\nabla g_j(y)\neq0$ whenever $y\in\K$ and $g_j(y)=0$.
If  Assumption \ref{geom} holds then $\K$ is convex and $\Omega$ in (\ref{thmain-2}) with $d:=\max_j d_j$,
is a semidefinite representation of $\K$.
\end{theorem}
\begin{proof}
That $\K$ is convex follows from Lemma \ref{lemmaconvex}.
We next prove that the PP-BDR property defined in Lasserre \cite{lasserre-sdr} holds for $\K$.
Let $f\in\R[X]$ be a linear polynomial with coefficient vector $\f\in\R^n$ (i.e., $X\mapsto f(X)=\f^TX$)
and consider the optimization problem $\P:\: \min \:\{\f^Tx\::\:x\in\K\}$.
As $\K$ is compact, let $x^*\in\K$ be a global minimizer of $f$. The Fritz-John optimality conditions state
that there exists $0\neq\lambda\in\R^{m+1}_+$ such that
\begin{equation}
\label{fj}
\lambda_0 \,\f=\sum_{j=1}^m\lambda_j \,\nabla g_j(x^*);\quad\lambda_j\,g_j(x^*)=0\quad\forall j=1,\ldots,m.\end{equation}
(See e.g. \cite{john}.)
We first prove by contradiction that if Slater's condition and the nondegeneracy assumption hold then 
$\lambda_0>0$. 
Suppose that $\lambda_0=0$ and let $J:=\{j\in\{1,\ldots,m\}\::\:\lambda_j>0\}$; hence $J$ is nonempty as $\lambda\neq0$. With $x_0\in\K$ such that $g_j(x_0)>0$ (as Slater's condition holds, one such $x_0$ exists),
let $B(x_0,\rho):=\{z \::\: \Vert z-x_0\Vert\leq\rho\}$. For $\rho$ sufficiently small,
$B(x_0,\rho)\subset\K$ and $g_j(z)>0$ for all $z\in B(x_0,\rho)$ and every $j=1,\ldots,m$. Then
by (\ref{fj}) and $\lambda_0=0$,
\[0=\sum_{j=1}^m\lambda_j \,\langle\nabla g_j(x^*),z-x^*\rangle,\qquad\forall z\in B(x_0,\rho),\]
which in turn implies (by nonnegativity of each term in the above sum) 
\[\langle\nabla g_j(x^*),z-x^*\rangle =0,\qquad \forall z\in B(x_0,\rho),\: j\in J.\]
%because $g_j(x^*)=0$ if $j\in J$ and (\ref{statconvex}) holds.
But this clearly implies $\nabla g_j(x^*)=0$ for every $j\in J$, in contradiction
with the nondegeneracy assumption. Hence $\lambda_0>0$ and by homogeneity,
we may and will take $\lambda_0=1$.

Therefore, letting $Y:=x^*$ in (\ref{everyj}), the polynomial $X\mapsto f(X)-f^*$ can be written
\begin{eqnarray*}
\f^TX-f^*&=&\sum_{j=1}^m\lambda_j\displaystyle\left[\:\langle \nabla g_j(x^*),X-x^*\rangle\:\right]\\
&=&\sum_{j=1}^m\lambda_j\left[\sum_{k=0}^m\sigma_{jk}(X,x^*)\, g_k(X)
+\sum_{k=0,k\neq j}^m\psi_{jk}(X,x^*) \,g_k(x^*)\right.\\
&&\left.+\,\psi_j(X,x^*)\,g_j(x^*)\right]
\end{eqnarray*}
where we have used (\ref{everyj}) with $Y=x^*$ and $\epsilon=0$.
Next, observe that :
\begin{eqnarray*}
X\mapsto \sigma_{jk}(X,x^*)&\in&\Sigma^2[X]\qquad\mbox{[as $\sigma_{jk}\in\Sigma^2[X,Y]$]}\\
X\mapsto \psi_{jk}(X,x^*) \,g_k(x^*)&\in&\Sigma^2[X]\qquad\mbox{[as $\psi_{jk}\in\Sigma^2[X,Y]$ and $g_j(x^*)\geq0$]}\\
\lambda_jg_j(x^*)&=&0\qquad j=1,\ldots,m.\end{eqnarray*}
And so, as $\lambda\in\R^m_+$,
\begin{equation}
\label{aux}
X\mapsto \f^TX-f^*\,=\,\Delta_0(X)+\sum_{j=1}^m\Delta_j(X)\,g_j(X),\end{equation}
for SOS polynomials $(\Delta_j)_{j=0}^m\subset\Sigma^2[X]$ defined by
\begin{eqnarray*}
X\mapsto\Delta_0(X)&=&\sum_{j=1}^m\lambda_j\left(\sum_{k=0,k\neq j}^m\psi_{jk}(X,x^*) \,g_k(x^*)\right)\\
X\mapsto\Delta_j(X)&=&\sum_{l=1}^m\lambda_l\,\sigma_{lj}(X,x^*),\qquad j=1,\ldots,m.
\end{eqnarray*}

%In addition, the degree of the SOS polynomials $(\Delta_j)\in\Sigma^2[X]$ is bounded by $2\max_jr_j$, independently of $\f$. 
Write every affine polynomial $f\in\R[X]$ as $\f^TX +f_0$ for some $\f\in\R^n$ 
and $f_0=f(0)$.
If $f$ is nonnegative on $\K$ then from (\ref{aux}),
\begin{eqnarray*}
f(X)\,=\,\f^TX-f^*+f^*+f_0&=&f^*+f_0+\Delta_0(X)+\sum_{j=1}^m\Delta_j(X)\,g_j(X)\\
&=&\widehat{\Delta}_0(X)+\sum_{j=1}^m\Delta_j(X)\,g_j(X)\qquad\forall X,
\end{eqnarray*}
with $\widehat{\Delta}_0\in\Sigma^2[X]$ (because $f^*+f_0\geq0$)
and so, the PP-BDR property holds for $\K$ with order $d$.
By \cite[Theor. 2]{lasserre-sdr}, $\K$ is SDr with the semidefinite representation (\ref{thmain-2}).
\end{proof}
\vspace{0.2cm}

We next show that the two sufficient conditions of strict convexity
and SOS-convexity of Helton and Nie \cite{HN1} 
in Theorem \ref{suffhn1} 
both imply that Assumption \ref{geom} holds and so Theorem \ref{thmain} contains
Theorem \ref{suffhn1} as  a special case.\\

\begin{corollary}
\label{finalsdr}
Let $\K$ in (\ref{setk}) be convex and both Assumption \ref{assput} and Slater's condition hold. Assume that either $-g_j$ 
is SOS-convex or $-g_j$ is convex on $\K$ and $-\nabla^2g_j\succ0$ on $\K\cap\{x\::\:g_j(x)=0\}$, for every $j=1,\ldots,m$. Then Assumption \ref{geom} holds
and so Theorem \ref{thmain} applies.
\end{corollary}
\begin{proof}
By Lemma \ref{prop2}, for every $j=1,\ldots,m$, write
\[(X,Y)\quad\mapsto \quad g_j(X)-g(Y)-\left\langle \nabla g_j(Y),X-Y\right\rangle\,=\,\]
\[\left\langle (X-Y),\underbrace{\left(\int_0^1\int_0^t\nabla^2
g_j(Y+s(X-Y))\,dsdt\right)}_{F_j(X,Y)}\,(X-Y)\right\rangle.\]
If $-\nabla^2g_j\succ0$ on $y\in\K$ with $g_j(y)=0$, then
from the proof of \cite[Lemma 19]{HN1}, $-F_j(x,y)\succ0$ for all
$x,y\in\K$ with $g_j(y)=0$. In other words,
%the compact set $\{x\in\K \,:\:g_j(x)=0\}$, 
$-F_j(x,y)\succeq\delta I_n$ on $\Omega_j$ (defined in (\ref{omegaj})) for some $\delta>0$.
Therefore, by the matrix polynomial version of Putinar Positivstellensatz
in \cite[Theor. 29]{HN1},
\begin{equation}
\label{hnaux}
-F_j(X,Y)\,=\,\sum_{k=0}^m\widehat{\sigma}_{jk}(X,Y)g_k(X)+
\sum_{k=0,k\neq j}^m\widehat{\psi}_{jk}(X,Y)g_k(Y)+
\widehat{\psi}_{j}(X,Y)g_j(Y)\end{equation}
for some SOS matrix polynomials
$(\widehat{\sigma}_{jk}(X,Y))$, $(\widehat{\psi}_{jk}(X,Y))$ and some 
matrix polynomial $\widehat{\psi}_{j}(X,Y)$.

On the other hand, if $-g_j$ is SOS-convex 
then by Lemma \ref{prop},
$-F_j(X,Y)$ is SOS and therefore (\ref{hnaux}) also holds (take 
$\widehat{\sigma}_{jk}\equiv 0$ for all $k\neq 0$,
$\widehat{\psi}_{jk}\equiv 0$ for all $k$ and $\widehat{\psi}_j\equiv 0$).
But then
\begin{eqnarray*}
g_j(X)-g(Y)-\left\langle \nabla g_j(Y),X-Y\right\rangle&=&
\left\langle (X-Y),F_j(X,Y)(X-Y)\right\rangle\\
&=&-\sum_{k=0}^m\left\langle (X-Y),\widehat{\sigma}_{jk}(X,Y)(X-Y)\right\rangle g_k(X)\\
&-&\sum_{k=0,k\neq j}^m\left\langle (X-Y),\widehat{\psi}_{jk}(X,Y)(X-Y)\right\rangle g_k(Y)\\
&-&\left\langle (X-Y),\widehat{\psi}_{j}(X,Y)(X-Y)\right\rangle g_j(Y)\\
&=&-\sum_{k=0}^m\sigma_{jk}(X,Y)\, g_k(X)-\\
&&\sum_{k=0,k\neq j}^m\psi_{jk}(X,Y) \,g_k(Y)
-\,\psi_j(X,Y)\,g_j(Y)
\end{eqnarray*}
for all $X,Y$ and for some SOS polynomials
$\sigma_{jk},\psi_{jk}\in\R[X,Y]$ and some polynomial $\psi_j\in\R[X,Y]$.
Equivalently,
\begin{eqnarray*}
\langle\nabla g_j(Y),X-Y)&=&g_j(X)-g_j(Y)
+\sum_{k=0}^m\sigma_{jk}(X,Y)\, g_k(X)\\
&&+\sum_{k=0,k\neq j}^m\psi_{jk}(X,Y) \,g_k(Y)
+\psi_j(X,Y)\,g_j(Y)\\
&=&\sum_{k=0}^m\sigma'_{jk}(X,Y)\, g_k(X)
+\sum_{k=0,k\neq j}^m\psi_{jk}(X,Y) \,g_k(Y)\\
&&+\psi'_j(X,Y)\,g_j(Y)
\end{eqnarray*}
for some SOS polynomials $\sigma'_{jk},\psi_{jk}\in\Sigma^2[X,Y]$
and some polynomial $\psi'_{j}\in\R[X,Y]$. 
In other words, Assumption \ref{geom} holds, which concludes the proof.
\end{proof}
\vspace{0.2cm}

Hence if each $-g_j$ is SOS-convex or convex on $\K$ with $-\nabla^2g_j\succ0 $ on $\K\cap\{x\::\:g_j(x)=0\}$,
one obtains a numerical scheme to obtain the parameter $d$ in Theorem
\ref{thmain} as well as the semidefinite representation 
(\ref{thmain-2}) of $\K$. Solve the semidefinite programs (\ref{test}) with degree
parameter $d_j$. Eventually, $\rho_j=0$ for every $j=1,\ldots,m$.\\

\begin{ex}
{\small{\rm Consider the convex set $\K$ in (\ref{setexample}) of Example \ref{ex1}
for which the defining polynomial $g_1$ of $\K$ is not concave.
We have seen that 
Assumption \ref{geom} holds (up to $\rho_1\approx 10^{-11}$, close to machine precision)
and $\max [d_1,d_2]=3$. By Theorem \ref{thmain}, if $\rho_1$ would be exactly $0$,
the set
\begin{equation}
\label{derniere}
\Omega\,:=\,\left\{(x,\y)\in\R^n\times \R^{s(6)}\::\:\left\{\begin{array}{ll}
M_{3}(\y)&\succeq0\\
M_{2}(g_j\,\y)&\geq 0,\quad j=1,2\\
L_\y(X_i)&=x_i,\quad i=1,2\\
y_0&=1\end{array}\right.\right..\end{equation}
would be a semidefinite representation of $\K$.

At least in practice, for every linear polynomial $f\in\R[X]$, minimizing
$L_\y(f)$ over $\Omega$
%under the constraints $M_3(\y),M_2(g_j\y)\succeq0$, $j=1,2$, and $y_0=1$,
yields the desired optimal value $f^*:=\min_{x\in\K} f(x)$, up to 
$\rho_1\approx -10^{-11}$.
%(close to machine precision) where $\rho_1$ was obtained in solving (\ref{test}) for $j=1$.

Indeed, let $f\in\R[X]$ be $\f^TX$ for some vector $\f\in\R^n$.
In minimizing $f$ over $\K$, one has $\f=\lambda_1\nabla g_1(x^*)+\lambda_2\nabla g_2(x^*)$ for some
$\lambda\in\R^2_+$, some $x^*\in\K$ with $\lambda_ig_i(x^*)=0$, $i=1,2$, and
$f^*=\lambda_1 \langle\nabla g_1(x^*),x^*\rangle+\lambda_2 \langle\nabla g_2(x^*),x^*\rangle
=\min_{x\in\K}\f^Tx$.
Let $x$ be as in (\ref{derniere}), arbitrary. Then
\[\f^Tx-f^*\,=\,L_\y(f(X)-f^*)\,=\,\sum_{i=1}^2
\lambda_iL_\y(\langle \nabla g_i(x^*),X-x^*\rangle).\]
If $\lambda_1>0$ so that $g_1(x^*)=0$, use (\ref{aux}) to obtain
\[L_\y(\langle \nabla g_1(x^*),X-x^*\rangle)=
L_\y(\rho_1+\Delta_0(X)+\sum_{j=1}^2\Delta_j(X)g_j(X))\geq \rho_1,\]
because $L_\y(\Delta_0)\geq0$ follows 
from $M_3(\y)\succeq0$, and $L_\y(\Delta_jg_j)\geq0$, $j=1,2$, follows from
$M_2(g_1\y),M_2(g_2\y)\succeq0$. If $\lambda_2>0$ so that $g_2(x^*)=0$, then from (\ref{newtest})
\[L_\y(\langle \nabla g_2(x^*),X-x^*\rangle)\,=\,L_\y(g_2(X)-
\langle (X-x^*),\nabla^2g_2(x^*)(X-x^*)\rangle)\geq0,\]
because $L_\y(g_2)\geq0$ follows from $M_2(g_2\,\y)\succeq0$ 
whereas the second term is nonnegative as
$\langle (X-x^*),-\nabla^2g_2(x^*)(X-x^*)\rangle$ is SOS and $M_3(\y)\succeq0$.
Hence $\f^Tx-f^*\geq\rho_1$. On the other hand, from
$\K\subseteq\{x\::\:(x,y)\in\Omega\}$, one finally obtains the desired result
\[f^*+\rho_1\,\leq\min\:\{\f^Tx\::\,(x,y)\in\Omega\}\,\leq\, f^*.\]
}}

\end{ex}

\section{Conclusion}
As well-known, convexity is a highly desirable property in optimization.
We have shown that it also has important specific consequences 
in polynomial optimization. For instance,  for polynomial optimization problems with 
SOS-convex or strictly convex polynomial data,  
the basic SDP-relaxations of the moment approach  \cite{lasserre1} 
{\it recognizes} convexity and finite convergence occurs.
Similarly, the set $\K$ has a semidefinite representation, explicit in terms of the defining 
polynomials $(g_j$).

The class of SOS-convex polynomials introduced in Helton and Nie \cite{HN1} is particularly interesting
because the semidefinite constraint to handle in the semidefinite relaxation only
involves the Hankel-like moment matrix
which does {\it not} depend on the problem data! Hence one might envision
a dedicated SDP solver that would take into account this peculiarity as Hankel-like or Toeplitz-like
matrices enjoy very specific properties.
Moreover, if restricted to this class of polynomials, Jensen's inequality can be 
extended to linear functionals in the dual cone of SOS polynomials (hence not necessarily
probability measures).

Therefore, a topic of further research is to evaluate how {\it large} is the subclass of SOS-convex 
polynomials in the class of convex polynomials, and if possible, to also provide simple sufficient conditions for SOS-convexity.

\section*{Acknowledgements}
The author wishes to thank L. Tuncel and Y. Nesterov for helpful discussions 
on various characterizations of convex sets, and also two anonymous referees 
for several corrections  as well as  suggestions and remarks 
to improve a first version of this paper.

\end{document}